\documentclass{amsart}

\newcommand     {\é}    {\'e}
\newcommand \lbb[1] {\label{#1}}

\newcommand {\1}    {{\bf 1}}

\newcommand     {\CD}   {{\C[\partial]}}
\newcommand     {\D}    {{\partial}}

\newcommand     {\C}    {{\mathbb C}}
\newcommand     {\N}    {{\mathbb N}}
\newcommand     {\Z}    {{\mathbb Z}}
\newcommand     {\G}    {{\mathfrak g}}

\newcommand     {\ad}   {\mbox{ad\,}}

\newcommand     {\End}  {\mbox{End}}

\newcommand     {\gc}   {{\mbox{\it gc}}}

\newcommand     {\Vir}  {{\mbox{\it Vir}}}

\newcommand     {\Tor}  {{\mbox{\it Tor }}}

\newcommand     {\Chom} {{\mbox{\it Chom}}}

\newcommand     {\Conf} {{\mbox{\it Cur }}}

\newcommand     {\rk}   {\mbox{rk\,}}

\newcommand     {\Rad}  {{\mbox{\it Rad\,}}}

\newcommand     {\id}   {\mbox{id}}
\newcommand     {\open} {``}
\newcommand     {\close}{''}
\newcommand     {\conf} {{\mbox{\scriptsize {\it Lie}}}}

\newtheorem{thm}{Theorem}[section]

\newtheorem{lemma}{Lemma}[section]
\newtheorem{prop}{Proposition}[section]
\newtheorem{cor}{Corollary}[section]

\theoremstyle{definition}

\newtheorem{defn}{Definition}[section]

\newtheorem{ex}{Example}[section]

\newtheorem{rem}{Remark}[section]

\numberwithin{equation}{section}
\begin{document}

\title[Finite Vertex Algebras and Nilpotence]
{Finite Vertex Algebras and Nilpotence}
\author[A.~D'Andrea]{Alessandro D'Andrea}
\thanks{The author was partially supported by
PRIN ``Spazi di Moduli e Teoria di Lie'' fundings from MIUR and
project MRTN-CT 2003-505078 ``LieGrits'' of the European Union}
\address{Dipartimento di Matematica, Universit\`a degli Studi di
Roma ``La Sapienza'', Roma}
\email{dandrea@@mat.uniroma1.it}
\date{March 22, 2007}


\begin{abstract}
I show that simple finite vertex algebras are commutative, and that
the Lie conformal algebra structure underlying a reduced (= without
nilpotent elements) finite vertex algebra is nilpotent.
\end{abstract}

\maketitle

\tableofcontents
\section{Introduction}

In this paper I investigate the effect of a finiteness assumption on
the singular part of the Operator Product Expansion of quantum
fields belonging to a vertex algebra. The vertex algebra structure
encodes algebraic properties of chiral fields in a 2-dimensional
Conformal Field Theory. Its axiomatic definition was given by
Borcherds in \cite{B}, and amounts to associating with each element
$v$ of a vector space $V$ a {\em vertex operator}, or quantum field,
$$Y(v,z) \in (\End V)[[z,z^{-1}]],$$
satisfying singular generalisations of the commutativity and unit
axioms for the left multiplication operators in an associative
algebra. This structure captures all algebraic properties of
families of mutually local fields, containing the identity field,
acting on some vector space of {\em physical states} $V$.

The products encoding the algebraic properties of a vertex algebra
structure are typically written in terms of a formal expansion
\begin{equation}\lbb{ope}
Y(a,z) Y(b,w) = \sum_{j=0}^\infty \frac{Y(a_{(j)}b, w)}{(z-w)^{j+1}}
+ :Y(a,z)Y(b,w):
\end{equation}
of the composition of quantum fields, called the Operator Product
Expansion (OPE for short). The singular part $$\sum_{j=0}^\infty
\frac{Y(a_{(j)}b, w)}{(z-w)^{j+1}}$$ of the OPE \eqref{ope} only
depends on the commutation properties of fields $Y(\cdot, z)$, and
is in many ways remindful of a Lie algebra; the regular part
$:Y(a,z) Y(b,w):$ is called {\em normally ordered product} or Wick
product, and essentially depends on the action of quantum fields on
$V$. It is to some extent similar to an associative commutative
product.

One can disregard the latter part of the structure, and axiomatize
the singular part of the OPE only. The thus obtained algebraic
structure is called \open (Lie) conformal algebra\close, but it is
also known in the literature as vertex Lie algebra \cite{P, DLM}, or
Lie pseudoalgebra \cite{BDK} over $\CD$. Lie conformal algebras were
introduced by Kac \cite{K} to characterize algebraic properties of
pairwise local formal distributions (in $z$ and $z^{-1}$) with
values in a Lie algebra. Hence every vertex algebra is thus also a
Lie conformal algebra, and the latter structure measures the failure
of the normally ordered product from being associative and
commutative.

Indeed, if the Lie conformal algebra underlying a vertex algebra is
trivial, then the product defined as $a \circ b = Y(a,z) b|_{z=0}$
gives \cite{B} a commutative associative algebra structure on $V$
from which it is possible to recover the vertex algebra product. In
fact, in this case $Y(a,z)b=(e^{zT}a)\circ b$ -- where $T$ is the
(infinitesimal) translation operator $Y(Ta, z) = dY(a,z)/dz$ -- and
is therefore completely determined by $\circ$.

The main interest in the role of Lie conformal algebras in vertex
algebra theory is due to the existence \cite{L, R, P, K} of a {\em
universal enveloping vertex algebra} functor, which is adjoint to
the forgetful functor from vertex to Lie conformal algebras. Many
interesting vertex algebras are in fact obtained as simple quotients
of the universal enveloping vertex algebra associated with a
suitably chosen (and typically much smaller) Lie conformal algebra.
However, a different strategy is possible: one could, in principle,
analyse possible vertex algebra structures by first studying the
singular OPE (i.e., the underlying Lie conformal algebra), and then
inserting the normally ordered product on top of this. My aim in
this note is to show that this strategy gives interesting results
when $V$ is a finitely generated $\C[T]$-module.

The physically interesting vertex algebras which are usually
considered are called Vertex Operator Algebras (=VOAs). They are
graded vector fields that are often endowed with some additional
structure (e.g., a Virasoro field inducing the grading); however,
all known examples are very large objects -- typically of
superpolynomial growth. A first explanation for this is Borcherds'
observation \cite{B2} that in a finite-dimensional vertex algebra
all products $Y(a,z)b$ are necessarily regular in $z$, as the
underlying Lie conformal algebra must be trivial (because all
elements are torsion, see \cite{K,DK}). The vertex algebra structure
reduces, as mentioned above, to that of a unital finite-dimensional
commutative associative algebra with a derivation $T$. Also, in the
presence of a grading induced by a Virasoro field the dimension of
the homogeneous component of degree $n$ is at least the number of
partitions of $n$, which grows super-polynomially. However, no
interesting (e.g., simple) examples are known of a vertex algebra
structure on a graded vector space of polynomial growth, even in the
absence of the additional requirements for a VOA.

After finite-dimensional vector spaces, the next to easiest case is
that of finitely generated modules over $\C[T]$: they are vector
spaces of linear growth when a grading is given, and this is the
case I handle in this note. There seems to be no previously known
description of algebraic properties of such {\em finite vertex
algebras} that are not finite-dimensional.

Section \ref{generalities} contains a list of definitions and basic
results in the theory of Lie conformal algebras and vertex algebras.
I also exhibit the motivating Example \ref{finitevertex}, showing
that there exist finite vertex algebras that do not reduce to
associative commutative algebras. However, the vertex algebra
provided in this example is constructed by means of nilpotent
elements: as in the case of commutative algebras, such elements form
an ideal of the vertex algebra which I call {\em nilradical}. The
quotient of $V$ by its nilradical has no non-zero (strongly)
nilpotent elements.

Sections \ref{wick} and \ref{solvability} apply results from
\cite{D} to the case of finite vertex algebras in order to show that
finite simple vertex algebras are commutative (Theorem
\ref{simpleiscommutative}), hence that the Lie conformal algebra
structure underlying a finite vertex algebra is always solvable
(Theorem \ref{finiteissolvable}).

In Section \ref{adjointdecomposition}, I study how a finite vertex
algebra decomposes under the adjoint action of a Lie conformal
subalgebra. After describing the representation theory of finite
solvable Lie conformal algebras I show, in Theorem
\ref{vertexideal}, that the generalized weight submodule with
respect to any non-zero weight is an abelian ideal of the vertex
algebra structure. The presence of abelian ideals witnesses the
existence of nilpotent elements, therefore there can be no non-zero
weights in the absence of nilpotent elements. This strong algebraic
fact is the basis for the results presented in Section
\ref{coreresults}, and is proved by means of the identity
\eqref{genwickeq} introduced in Section \ref{wick}.

The main result from last section is Theorem \ref{key} stating that
any element $s$ lying in a finite vertex algebra $V$ with trivial
nilradical has a nilpotent adjoint conformal action on $V$. By a Lie
conformal algebra analogue of Engel's theorem, developed in
\cite{DK}, the Lie conformal algebra underlying $V$ must indeed
(Theorem \ref{main}) be nilpotent. This statement essentially
depends on both the finiteness assumption and the presence of a
vertex algebra structure: it basically means that every finite
vertex algebra may be described as an extension of a nilpotent (as a
Lie conformal algebra) vertex algebra by an ideal only containing
nilpotent elements (i.e., contained in the nilradical).

The representation theory of the Virasoro Lie algebra or of affine
Kac-Moody algebras is never used. The spirit of this paper is that
the interplay between the Operator Product Expansion and the
$\lambda$-bracket, in the case of a vertex algebra linearly
generated by a finite number of quantum fields together with their
derivatives, is strong enough to allow one to prove a number of
results in a totally elementary way, even in the absence of a
grading.

Some of the ideas contained in this work originate from an old
manuscript written between 1998 and 1999 while I was visiting
Universit\é de Paris VI \open Pierre et Marie Curie\close\, and
Universit\é de Strasbourg \open Louis Pasteur\close\, as a European
Union TMR post-doc. I would like to thank both institutions for
hospitality.

\section{Generalities on vertex and Lie conformal
algebras}\lbb{generalities}
\subsection{Vertex algebras}

In what follows I quote some well-known facts about vertex algebras:
precise statements and proofs can be found in \cite{K}. Let $V$ be a
complex vector space. A {\em field} on $V$ is an element $\phi(z)
\in (\End V)[[z, z^{-1}]]$ with the property that $\phi(v) \in
V((z)) = V[[z]][z^{-1}]$ for every $v \in V$. In other words, if
$$\phi(z) = \sum_{i \in \Z} \phi_i z^{-i-1}$$ then $\phi_n(v) = 0 $
for sufficiently large $n$.

\begin{defn}
A vertex algebra is a (complex) vector space $V$ endowed with a
linear {\em state-field correspondence} $Y:V \to (\End V)[[z,
z^{-1}]]$, a {\em vacuum element} $\1$ and a linear endomorphism
$T \in \End V$ satisfying the following properties:
\begin{itemize}
\item {\bf Field axiom}: $Y(v,z)$ is a field for all $v\in V$.
\item {\bf Locality axiom}: For every $a, b \in V$
$$(z-w)^N[Y(a,z), Y(b,w)] = 0$$ for sufficiently large $N$.
\item {\bf Vacuum axiom}: The vacuum element $\1$ is such that
$$Y(\1,z) = \id_V,\qquad Y(a,z)\1 \equiv a \mod
zV[[z]],$$ for all $a \in V$.
\item {\bf Translation invariance}: $T$ satisfies $$[T, Y(a,z)]
= Y(Ta, z) = \frac{d}{dz}Y(a,z),$$ for all $a\in V$.
\end{itemize}
\end{defn}
Note that the vector space $V$ carries a natural $\C[T]$-module
structure. Fields $Y(a,z)$ are called {\em vertex operators}, or
{\em quantum fields}.

A vertex algebra is a family of pairwise local fields acting on $V$
containing the identity (constant) field. Indeed, every family of
pairwise local fields containing the identity field can be realized
as a vertex algebra up to changing the vector space $V$ of physical
states (see \cite{K}). The vertex algebra structure therefore
captures all algebraic aspects of families of pairwise local fields.
A vertex algebra $V$ is {\em finite} if $V$ is a finitely generated
$\C[T]$-module.

There are two basic constructions of new vertex operators from two
given ones. The first one is given by rephrasing what we earlier
called \open singular OPE\close: since $(z-w)^N$ kills the
commutator $[Y(a,z), Y(b,w)]$, the latter may be expanded into a
linear combination:
\begin{equation}
\sum_{j=0}^{N-1} c_j(w) \frac{\delta^{(j)} (z-w)}{j!},
\end{equation}
where $$\delta(z-w)= \sum_{j \in \Z} w^j z^{-j-1}$$ is the Dirac
delta formal distribution and $\delta^{(j)}$ its $j$.th derivative
with respect to $w$. The uniquely determined fields $c_j(w)$ are
then vertex operators $Y(c_j,w)$ corresponding to elements $c_j =
a_{(j)} b= a_{(j)}(b)$ where the $a_{(j)}\in \End V$ are the
coefficients of
$$Y(a,z) = \sum_{j\in \Z} a_{(j)} z^{-j-1}.$$
It is customary to view the $\C$-bilinear maps $a \otimes b \mapsto
a_{(j)} b, j \in \Z,$ as products describing the vertex algebra
structure. Locality can be rephrased by stating that commutators
between coefficients of quantum fields satisfy:
\begin{equation}\lbb{liebracket}
[a_{(m)}, b_{(n)}] = \sum_{j \geq 0} {m \choose j}
(a_{(j)}\,b)_{(m+n-j)},
\end{equation}
for all $a, b \in V$.

Another way to put together quantum fields to produce new ones is
given by the normally ordered product (or Wick product) defined as:
\begin{equation}
:Y(a,z)Y(b,z): \,\,= Y(a,z)_+Y(b,z)+ Y(b,z)Y(a,z)_-,
\end{equation}
where
\begin{equation}
Y(a,z)_- = \sum_{j\in\N} a_{(j)} z^{-j-1},\,\,\,\,\, Y(a,z)_+ =
Y(a,z)-Y(a,z)_-.
\end{equation}
Then $:Y(a,z)Y(b,z):$ is also a vertex operator, and it equals
$Y(a_{(-1)}b,z)$.

\begin{ex}\lbb{holomorphic}
Let $V$ be a unital associative commutative algebra, $T$ a
derivation of $V$. Then setting $Y(a,z)b = (e^{zT}a)b$ and
choosing the unit $\1 \in V$ to be the vacuum element makes $V$
into a vertex algebra.
\end{ex}
Such a vertex algebra is called {\em holomorphic}\footnote{Notation
and terminology are often contrasting in the vertex algebra world,
and this is no exception. Notice that in most of the literature
vertex operator algebras are known to be holomorphic if they have a
semi-simple representation theory, and the adjoint representation is
the unique irreducible module, see \cite{blabla}.} in \cite{K}, and
is the \open uninteresting\close\, case of a vertex algebra
structure. It occurs whenever all vertex operators are regular in
$z$: in this case one can always construct an associative
commutative algebra, together with a derivation, inducing the vertex
algebra structure as in Example \ref{holomorphic}. This is always
the case when $V$ is finite-dimensional \cite{B, B2}.
One of the consequences of the vertex algebra axioms is the
following:
\begin{itemize}
\item {\bf Skew-commutativity}: \qquad $Y(a,z)b = e^{zT} Y(b,-z)a$
\end{itemize}
for all choices of $a,b$.

If $A$ and $B$ are subsets of $V$, then we may define $A \cdot B$ as
the $\C$-linear span of all products $a_{(j)} b$, where $a\in A, b
\in B, j \in \Z$. If $B$ is a $\C[T]$-submodule of $V$, then $A\cdot
B$ is also a $\C[T]$-submodule of $V$, as by translation invariance
$T$ is a derivation of all $j$-products, and $(Ta)_{(j)} = - j
a_{(j-1)}$.

Notice that by skew-commutativity, $A \cdot B$ is contained in the
$\C[T]$-sub\-module generated by $B \cdot A$; equality $A \cdot B =
B \cdot A$ then holds whenever $A \cdot B$ and $B \cdot A$ are both
$\C[T]$-submodules of $V$. Also, observe that $A \subset A \cdot V$
by the vacuum axiom and that $A \cdot V$ is always a
$\C[T]$-submodule of $V$, as $a_{(-2)} \1 = T a$. In particular, $a
\cdot V = \C a \cdot V$ is a $\C[T]$-module of $V$ containing $a$.

Before proceeding, recall that a {\em subalgebra} of a vertex
algebra $V$ is a $\C[T]$-submodule $U$ containing $\1$ such that $U
\cdot U = U$. In other words, all coefficients of $Y(a,z)b$ belong
to $U$ whenever $a$ and $b$ do. Similarly, a $C[T]$-submodule $I$ of
$V$ is an {\em ideal} if $I \cdot V \subset I$; skew-commutativity
then shows that $V \cdot I = I \cdot V$. A proper ideal can never
contain the vacuum $\1$, and if $M$ is an ideal of $V$, then $M +
\C\1$ is a subalgebra, whose rank as a $\C[T]$-module equals that of
$M$.
%

The quotient $V/I$ of a vertex algebra $V$ by an ideal $I$ has a
unique vertex algebra structure making the canonical projection
$\pi: V \to V/I$ a vertex algebra homomorphism, i.e., a
$\C[T]$-homomorphism such that $\pi(a_{(n)} b) = \pi(a)_{(n)}
\pi(b)$ for every $a, b\in V$, $n\in \Z$.

A vertex algebra $V$ is {\em commutative} if all quantum fields
$Y(a, z), a \in V$ commute with one another; equivalently, if
$a_{(n)}b = 0$ for all $a, b \in V, n \geq 0$. Commutative vertex
algebras are all as in Example \ref{holomorphic}. The {\em centre}
of $V$ is the subspace of all elements $c\in V$ such that $a_{(n)}c
= 0 = c_{(n)}a$ for all $a \in V, n \geq 0$. Then, by
\eqref{liebracket}, coefficients of $Y(c, z)$ commute with
coefficients of all quantum fields.

\begin{lemma}\lbb{kercisideal}
Assume that $c$ lies in the centre of a vertex algebra $V$. Then the
subspace $\ker c_{(-1)}$ is stable under the action of coefficients
of all quantum fields $Y(a, z), a \in V$. In particular, $\ker
c_{(-1)}$ is an ideal of $V$ as soon as it is a $\C[T]$-submodule,
e.g., when $T c = 0$.
\end{lemma}
\begin{proof}
By \eqref{liebracket} we have
$$a_{(m)}(c_{(-1)}x) - c_{(-1)}(a_{(m)}x) = \sum_{j \geq 0} {m \choose j}
(a_{(j)}c)_{(m-j-1)}x,$$ for all $a \in V, m \in \Z$. Since $c$ lies
in the centre of $V$, the right hand side vanishes; hence, if $x \in
\ker c_{(-1)}$, it follows that $a_{(m)}x \in \ker c_{(-1)}$ as
well. The last claim follows from the fact that $c_{-1}(T x) = T
(c_{(-1)}x) - (T c)_{(-1)}x$.
\end{proof}


\subsection{A non-commutative finite vertex algebra}

The following is an example of a vertex algebra structure on a
finitely generated $\C[T]$-module for which some positive products
$u_{(j)} v, j \geq 0$ are non-zero.
\begin{ex}\lbb{finitevertex}
Let $V = \C[T]a \oplus \C[T]b \oplus \C \1$. Define $\1$ to be the
vacuum element of $V$ and set $$Y(a,z)b = Y(b,z)a = Y(b,z)b = 0,$$
$$Y(\1, z) = \id_V$$
$$Y(a,z)\1 = e^{zT} a \qquad Y(b,z)\1 = e^{zT} b$$
$$Y(a,z)a = e^{zT/2} \psi(z) b,$$ where $\psi(z) = \psi(-z)$ is
any Laurent series in $z$. Extend by $\C$-linearity the state-field
correspondence $Y$ to all of $V$ after setting:
$$Y(Tu,z)v = \frac{d}{dz}Y(u,z)v$$
and
$$Y(u,z)(Tv) = (T - \frac{d}{dz})(Y(u,z)v),$$
so that translation invariance is satisfied. The only vertex algebra
axiom still to check is locality, and the only non-trivial statement
to prove is
$$(z-w)^n [Y(a,z), Y(a,w)]\1 = 0,$$
for some $n$. However we have
\begin{align*}
[Y(a,z), Y(a,w)]\1 & = Y(a,z) e^{wT}a - Y(a,w)e^{zT}a\\
& = e^{(z+w)T/2} \left(\iota_{z,w}\psi(z-w) -
\iota_{w,z}\psi(w-z)\right)b,
\end{align*}
where $\iota_{z,w}$ (resp. $\iota_{w,z}$) is a prescription to
consider the expansion in the domain $|z|>|w|$ (resp. in the domain
$|w|>|z|$), see \cite{K}. If we choose $n$ so that $z^n \psi(z)$ is
regular in $z$, multiplication by $(z-w)^n$ makes the above
expression zero, due to the fact that $\psi(z) = \psi(-z)$: it is an
{\em expansion of zero} in the sense of \cite{FHL}. Notice that if
we choose a non-regular $\psi(z)$, then at least one of the products
$a_{(j)}a, j\geq 0,$ is non-zero, so that $V$ is non-commutative.
\end{ex}

In the example above, elements $a$ and $b$ are {\em nilpotent}, in a
sense that we are about to clarify.

\subsection{The nilradical}

An ideal $I$ of a vertex algebra $V$ is {\em abelian} if $I^2 = I
\cdot I = 0$. An element $a \in V$ is {\em strongly nilpotent} of
degree $n$ if every product of elements in $V$ containing $a$ at
least $n$ times, under all $j$-products and any parenthesization,
gives $0$.

Let $x \in V$ be a strongly nilpotent element of degree $n>2$, and
$a$ be a non-zero product of $[(n+1)/2]$ copies of $x$. Then $a$ is
strongly nilpotent of degree two.

\begin{lemma}\lbb{squaretozero}
Let $a \in V$ be a strongly nilpotent element of degree two. Then
$a$ generates an abelian ideal of $V$.
\end{lemma}
\begin{proof}
Clear.
\end{proof}

\begin{rem}
An element $a\in V$ is (non-strongly) nilpotent of degree $n$ if
every product of at least $n$ copies of $a$, under any product and
parenthesization gives $0$. If $V$ is either commutative or graded,
then it is easy to show that every nilpotent element is strongly
nilpotent; this holds in general, as can be proved using
\cite[Remark 7.8]{BK}. In particular, $Y(a, z)a = 0$ guarantees that
$a \cdot V$ is an abelian ideal of $V$.
\end{rem}

%
%

\begin{cor}
A vertex algebra $V$ possesses non-zero strongly nilpotent elements
if and only if it contains a non-zero abelian ideal.
\end{cor}
\begin{proof}
Every non-zero element in an abelian ideal is strongly nilpotent of
degree $2$. The converse is Lemma \ref{squaretozero}.
\end{proof}
Let us now denote
$$I^1 = I, \qquad I^{n+1} = I^n \cdot I^n,\,\, n > 0.$$
Then $I$ is a {\em nil-ideal} if $I^n = 0$ for sufficiently large
values of $n$.
\begin{lemma}
Let $V$ be a vertex algebra, $N \subset V$ a nil-ideal, $\pi: V \to
V/N$ the natural projection. Then $I \subset V$ is a nil-ideal if
and only if $\pi(I)$ is a nil-ideal of $V/N$.
\end{lemma}
\begin{proof}
We have $\pi(I^n) = \pi(I)^n$, hence if $I$ is a nil-ideal, $\pi(I)$
is too. On the other hand, if $\pi(I)^n = \pi(I^n) = 0$, then $I^n
\subset N$. Thus, $I^{n+k} \subset N^k$ which is $0$ for
sufficiently large $k$.
\end{proof}
\begin{cor}\lbb{sumofnil}
The sum of nil-ideals is a nil-ideal.
\end{cor}
\begin{proof}
Let $I, J$ be nil-ideals of $V$, and let $\pi: V \to V/I$ be the
natural projection. Then $\pi(I + J)$ equals $\pi(J)$ which is a
nil-ideal.
\end{proof}
\begin{cor}
Let $V$ be a finite vertex algebra. Then $V$ has a unique maximal
nil-ideal.
\end{cor}
\begin{proof}
Existence of some maximal nil-ideal follows from finiteness of the
$\C[T]$-module $V$, which is therefore Noetherian; uniqueness from
Corollary \ref{sumofnil}.
\end{proof}
The unique maximal nil-ideal of a finite vertex algebra $V$ is
called the {\em nilradical} $N(V)$ of $V$. It is clear that every
strongly nilpotent element of $V$ lies in $N(V)$. Furthermore, the
quotient $V/N(V)$ has no strongly nilpotent elements, hence it has a
trivial nilradical. We will call a vertex algebra with a trivial
nilradical, or equivalently with no non-trivial abelian ideal, a
{\em reduced} vertex algebra.
\begin{rem}
It is true in general that an element lies in the nilradical $N(V)$
if and only if it is nilpotent, hence strongly nilpotent; however,
we will not need this fact.
\end{rem}

\subsection{Lie conformal algebras}\lbb{confalg}

Algebraic properties of commutators of quantum fields are encoded in
the notion of Lie conformal algebra. See \cite{DK} for generalities
on (Lie) conformal algebras and $\lambda$-brackets.

\begin{defn} A {\em Lie conformal algebra} is a
$\C[\partial]$-module $R$ with a $\C$-bilinear product $(a,b)\mapsto
[a\,_\lambda\, b] \in V[\lambda]$ satisfying the following axioms:
\begin{enumerate}
\item[{\bf (C1)}] $[a \,_\lambda\, b] \in R[\lambda],$

\item[{\bf (C2)}] $[\partial a \,_\lambda\, b] = -\lambda [a \,_\lambda\,
b],\,\,\, [a \,_\lambda\, \partial b] = (\partial +\lambda)
[a\,_\lambda\,b],$

\item[{\bf (C3)}] $[a \,_\lambda\, b] = - [b \,_{-\partial -\lambda}\, a],$

\item[{\bf (C4)}] $[a \,_\lambda\, [b \,_\mu\, c]] - [b \,_\mu\,
[a\,_\lambda\,c]] = [[a \,_\lambda\, b] \,_{\lambda + \mu}\, c],$
\end{enumerate}
\noindent for every $a,b,c\in V$.
\end{defn}

Any vertex algebra $V$ can be given a $\CD$-module structure by
setting $\D = T$. Then defining
$$ [a\,_\lambda\, b] = \sum_{n\in \N} \frac{\lambda^n}{n!}
a_{(n)}b$$ endows $V$ with a Lie conformal algebra structure.
Indeed, (C1) follows from the field axiom, (C2) from translation
invariance, (C3) from skew-commutativity, and (C4) from
\eqref{liebracket}. In all that follows we will denote the
infinitesimal translation operator $T$ in a vertex algebra by $\D$.

If $A$ and $B$ are subspaces of a Lie conformal algebra $R$, then we
may define $[A ,B]$ as the $\C$-linear span of all
$\lambda$-coefficients in the products $[a\,_\lambda b]$, where
$a\in A, b \in B$. It follows from axiom (C2) that if $B$ is a
$\CD$-submodule of $R$, then $[A, B]$ is also a $\CD$-submodule of
$R$. Notice that if $A$ and $B$ are both $\CD$-submodules, then $[A,
B] = [B, A]$ by axiom (C3). A {\em subalgebra} of a Lie conformal
algebra $R$ is a $\CD$-submodule $S \subset R$ such that $[S, S]
\subset S$.

A Lie conformal algebra $R$ is {\em solvable} if, after defining
$$R^{(0)} = R,\qquad R^{(n+1)} = [R^{(n)}, R^{(n)}], n\geq 0,$$
we find that $R^{(N)} = 0$ for sufficiently large $N$. $R$ is
solvable iff it contains a solvable ideal $S$ such that $R/S$ is
again solvable. Solvability of a nonzero Lie conformal algebra $R$
trivially fails if $R$ equals its {\em derived subalgebra} $R' = [R,
R]$. Similarly, $R$ is {\em nilpotent} if, after defining
\begin{equation}\lbb{centralseries}
R^{[0]} = R,\qquad R^{[n+1]} = [R, R^{[n]}], n\geq 0,
\end{equation}
we find that $R^{[N]} = 0$ for sufficiently large $N$.

An {\em ideal} of a Lie conformal algebra $R$ is a $\CD$-submodule
$I\subset R$ such that $[R, I] \subset I$. If $I, J$ are ideals of
$R$, then $[I, J]$ is an ideal as well. An ideal $I$ is said to be
{\em central} if $[R, I] = 0$, i.e., if it is contained in the {\em
centre} $Z=\{r \in R| [r\,_\lambda s] = 0 \mbox{ for all } s \in
R\}$ of $R$. $R$ is {\em abelian} if it coincides with its centre,
i.e. if $[R, R] = 0$.

A Lie conformal algebra $R$ is {\em simple} if its only ideals are
trivial, and $R$ is not {\em abelian}. An interesting such example
occurs when, for each choice of $0 \neq r \in R$, it occurs that
$[r, R] = [\C r, R] = R$. In this case $R$ is a {\em strongly
simple} Lie conformal algebra.

Notice that, when $V$ is a vertex algebra, we should distinguish
between ideals of the vertex algebra structure and ideals of the
underlying Lie conformal algebra. Indeed, ideals of the vertex
algebra are also ideals of the Lie conformal algebra, but the
converse is generally false, as it can be seen by noticing that
$\C\1$ is always a central ideal of the Lie conformal algebra
structure, but it is never an ideal of the vertex algebra.

In order to avoid confusion, we will denote by $V^\conf$ the Lie
conformal algebra structure underlying a vertex algebra $V$;
similarly, if $S \subset V$ is a $\CD$-submodule closed under all
nonnegative products $\,_{(n)}, n \in \N$, we will denote by
$S^\conf$ the corresponding Lie conformal algebra structure. The
reader should pay special attention to the fact that a vertex
algebra $V$ is commutative if and only if the Lie conformal algebra
$V^\conf$ is abelian, and that claiming that $I$ is an abelian ideal
of $V$ is a stronger statement than saying that $I$ is an abelian
ideal in $V^\conf$. We will say that $V$ is solvable (resp.
nilpotent), whenever $V^\conf$ is.

\subsection{Finite simple Lie conformal algebras}

Every Lie conformal algebra $R$ has a maximal solvable ideal, called
{\em radical} of $R$ and denoted by $\Rad R$. A Lie conformal
algebra is called semi-simple if it has no solvable ideal; the
quotient $R/\Rad R$ is always semi-simple.

An investigation of Lie conformal algebra structures on finitely
generated $\CD$-modules was undertaken in \cite{DK}, where a
classification of simple and semi-simple ones, together with
generalizations of standard theorems in Lie representation theory,
are presented.

It turns out that the only (up to isomorphism) simple Lie conformal
algebra structures over finitely generated $\CD$-modules are the
Virasoro conformal algebra and current conformal algebras over a
finite-dimensional simple Lie algebra, which are described below.
Semi-simple instances are direct sums of Lie conformal algebras that
are either simple or non-trivial semi-direct sums of a Virasoro
conformal algebra with a simple current one.

\begin{ex}
Let $R$ be a free $\CD$-module of rank one, generated by an element
$L$. Then
\begin{equation}\lbb{virasoro}
[L_\lambda L] = (\partial+2 \lambda) L
\end{equation}
uniquely extends to a Lie conformal algebra structure on $R$, which
is easily seen to be strongly simple. $R=\Vir$ is called {\em
Virasoro conformal algebra}.
\end{ex}

\begin{ex}
Let $\G$ be a finite-dimensional complex Lie algebra, and let $R =
\CD \otimes \G$. There exists a unique $\lambda$-bracket on $R$
extending
\begin{equation}\lbb{current}
[g_\lambda h] = [g,h],
\end{equation}
for $g,h \in \G \simeq 1 \otimes \G \subset R$, and satisfying all
axioms for a Lie conformal algebra. $R$ is called {\em current
conformal algebra} and is denote by $\Conf \G$. It is a simple Lie
conformal algebra whenever $\G$ is a simple Lie algebra. However,
$\Conf \G$ is never strongly simple, as for no choice of $g \in \G$
does $\ad g \in \End \G$ satisfy surjectivity.
\end{ex}

\subsection{Centre and torsion}

We now need a statement on subalgebras of a finite Lie conformal
algebra $R$ which have the same rank as $R$.


\begin{defn}
Let $U, V$ be $\CD$-modules. A {\em conformal linear map} from $U$
to $V$ is a $\C$-linear map $f_\lambda :U \to V[\lambda]$ such that
$f_\lambda (\D u) = (\D+\lambda) f_\lambda u$ for all $u \in U$.
\end{defn}
The space of all conformal linear maps from $U$ to $V$ is denoted by
$\Chom(U, V)$. It can be turned into a $\CD$-module via $$(\D
f)_\lambda u = -\lambda f_\lambda u.$$

\begin{rem}
Let $U, V, W$ be $\CD$-modules. A $\CD$-linear homomorphism $\phi: V
\to W$ induces a corresponding $\CD$-homomorphism $\phi_* :
V[\lambda] \to W[\lambda]$. Then if $f \in \Chom(U, V)$, the
composition $\phi_* \circ f$ lies in $\Chom(U, W)$.
\end{rem}

\begin{lemma}\lbb{quotient}
Let $U, V$ be $\CD$-modules, $f \in \Chom(U, V)$. If $U_0 \subset U,
V_0 \subset V$ are $\CD$-submodules such that $f_\lambda (u_0) \in
V_0[\lambda]$ for all $u_0 \in U_0$, then $f$ induces a unique
$\overline f \in \Chom(U/U_0, V/V_0)$.
\end{lemma}
\begin{proof}
Let $\pi: V \to V/V_0$ be the natural projection. Then $\pi_* \circ
f$ is a conformal linear map from $U$ to $V/V_0$ which kills all
elements from $U_0$.
\end{proof}

The most typical example of a conformal linear map comes from the
adjoint action in Lie conformal algebras. Indeed, if $R$ is a Lie
conformal algebra, and $r\in R$, then
$$(\ad r)_\lambda x = [r\,_\lambda x]$$
defines a conformal linear map from $R$ into itself, due to axiom
(C2).

\begin{lemma}[\cite{DK}]\lbb{torsion}
If $f \in \Chom(U, V)$ and $u \in \Tor U$, then $f_\lambda u = 0$.
\end{lemma}

\begin{cor}
The torsion of a Lie conformal algebra is contained in its centre.
\end{cor}
\begin{proof}
Let $R$ be a Lie conformal algebra, $r \in R, t \in \Tor R$. The
adjoint action of $r$ is a conformal linear map from $R$ into
itself, hence it maps the torsion element $t$ to $[r_\lambda t] = 0$
by Lemma \ref{torsion}.
\end{proof}

\begin{lemma}\lbb{cotorsion}
Let $S \subset R$ be finite Lie conformal algebras, such that $R/S$
is a torsion $\CD$-module. Then $S$ is an ideal of $R$ containing
$R'$, i.e., $R/S$ is abelian.
\end{lemma}
\begin{proof}
Since $S$ is a subalgebra of $R$, the adjoint action of $S$ on $R$
stabilizes the $\CD$-submodule $S$. By Lemma \ref{quotient}, $S$
acts on the quotient $R/S$, which is torsion. By Lemma
\ref{torsion}, the action of $S$ on $R/S$ is trivial, or in other
words $[S, R] \subset S$, which amounts to saying that $S$ is an
ideal of $R$.

Thus, the adjoint action of $R$ on itself stabilizes $S$, and we may
repeat the above argument to conclude that $R' = [R, R] \subset S$.
It immediately follows that $R/S$ is abelian.
\end{proof}

\subsection{Irreducible central extensions of the Virasoro conformal algebra}

In this section we compute all finite irreducible central extensions
of the Virasoro conformal algebra. Finite central extensions of
$\Vir$ are described, up to equivalence, by cohomology classes
\cite{BKV} of $H^2(\Vir, Z)$, where $Z$ is the finitely generated
$\CD$-module describing the centre\footnote{It is clear that every
finite central extension of $\Vir$ splits as an extension of
$\CD$-modules, as $\Vir$ is a free $\CD$-module.}.

The centre $Z$ being a finitely generated $\CD$-module, we can
decompose it (non-canonically) into a direct sum of its torsion with
a free $\CD$-module. This leads to a corresponding direct sum
decomposition of the related cohomology. In order to understand
finite central extensions of $\Vir$, it is thus sufficient to
compute $H^2(\Vir, Z)$ when $Z$ is either a free $\CD$-module of
rank one, or an indecomposable torsion $\CD$-module. The following
facts were proved in \cite{DK} and \cite{BKV} respectively:
\begin{prop}
All central extensions of $\Vir$ by a free $\CD$-module of rank one
are trivial.
\end{prop}
\begin{prop}\lbb{onedimcentralextensions}
Let $\C_\alpha, \alpha \in \C,$ be the $1$-dimensional $\CD$-module
on which the action of $\D$ is given via scalar multiplication by
$\alpha$. Then:
\begin{itemize}
\item if $\alpha \neq 0$ all central extensions of $\Vir$ by
$\C_\alpha$ are trivial;
\item if $\alpha = 0$ then there is a unique (up to isomorphism and
scalar multiplication) non-trivial central extension of $\Vir$ by
$\C = \C_0$ given by
\begin{equation}\lbb{irreducible}
[L_\lambda L] = (\D + 2 \lambda) L + \lambda^3.
\end{equation}
\end{itemize}
\end{prop}
\begin{rem}\lbb{virc0}
A computation of $2$-cocycles of $\Vir$ with values in the trivial
$\Vir$-module $\C_0$ shows that they are of the form $p(\lambda) =
c_1 \lambda + c_3 \lambda^3$, whereas trivial $2$-cocycles (i.e.,
$2$-coboundaries) are of the form $p(\lambda) = c_1 \lambda$.
\end{rem}

Recall that a central extension is called {\em irreducible} if it
equals its derived algebra. Clearly, no non-zero trivial central
extension is irreducible. My aim is to show that the non-trivial
central extension \eqref{irreducible} is the unique (non-zero)
irreducible finite central extension of $\Vir$.

\begin{prop}\lbb{trivial}
Let $C$ be a finitely generated torsion $\CD$-module on which $\D$
acts invertibly. Then every central extension of $\Vir$ by $C$ is
trivial.
\end{prop}
\begin{proof}
Let the central extension be given by
$$[L_\lambda L] = (\D + 2 \lambda)L + p(\lambda),$$
for some $p(\lambda) \in C[\lambda]$. By a computation similar to
that in \cite[Lemma 8.11]{DK}, one obtains $\D p(\lambda) = (\D +
2\lambda) p(0)$, whence $p(\lambda) = (\D + 2\lambda) \D^{-1} p(0)$.
Then $L + \D^{-1} p(0)$ is a standard generator of a Virasoro
conformal algebra, hence it splits the central extension.
\end{proof}

\begin{lemma}\lbb{centralextensions}
Solutions $p(\D, x) \in \C[\D, x]/(\D^{N+1})$ of
\begin{equation}\lbb{twococycle}
(\lambda - \mu) p(\D, \lambda + \mu) = (\D + \lambda + 2 \mu) p(\D,
\lambda) - (\D + 2 \lambda + \mu) p(\D, \mu) \mod \D^{N+1}
\end{equation}
are all of the form $p(\D, \lambda) = (\D + 2\lambda) q(\D) + c
\lambda^3 \D^N \mod \D^{N+1}, c \in \C$.
\end{lemma}
\begin{proof}
By induction on $N\geq 0$. The basis of induction follows from
Remark \ref{virc0}. Next assume $N>0$. Then \eqref{twococycle} also
holds modulo $\D^N$, and inductive assumption gives $p(\D, \lambda)
= (\D + 2 \lambda) q(\D) + c_0 \lambda^3 \D^{N-1} \mod \D^N$. As a
consequence:
$$p(\D, \lambda) = (\D + 2 \lambda) q(\D) + c_0 \lambda^3 \D^{N-1} +
\alpha(\lambda) \D^N \mod \D^{N+1}.$$ We can substitute this into
\eqref{twococycle} and get
$$(\lambda - \mu)\alpha(\lambda + \mu) - (\lambda + 2
\mu)\alpha(\lambda) + (2\lambda + \mu)\alpha(\mu) = c_0 (\lambda^3 -
\mu^3).$$ The left-hand side is linear in $\alpha$ and homogeneous
with respect to the joint degree in $\lambda$ and $\mu$. Hence we
can solve it degree by degree, looking for solutions of the form
$\alpha(x) = ax^n$. It is then easy to check that solutions only
exist when $c_0 = 0$, and are of the form $\alpha(x) = q_N \lambda +
c \lambda^3$. We conclude that
\begin{align*}
p(\D, \lambda) & = (\D + 2 \lambda) q(\D) + (q_N \lambda + c
\lambda^3) \D^N\\
& = (\D + 2 \lambda) \left(q(\D) + \frac{q_N}{2} \D^N\right) + c
\lambda^3 \D^N \mod \D^{N+1}.
\end{align*}
\end{proof}

\begin{prop}\lbb{jordanblock}
Let $C_N$ be a finitely generated torsion $\CD$-module isomorphic to
$\CD/(\D^{N+1})$, $N \geq 1$. Then there is a unique (up to
isomorphism and scalar multiplication) non-trivial central extension
of $\Vir$ by $C_N$ given by $[L_\lambda L] = (\D + 2 \lambda) L +
\lambda^3 \D^N.$
\end{prop}
\begin{proof}
The $2$-cocycle property for $p(\D, \lambda)$ as in
$$[L_\lambda L] =
(\D + 2 \lambda) L + p(\D, \lambda),$$ leads to solving
\eqref{twococycle}; hence Lemma \ref{centralextensions} gives $p(\D,
\lambda) = (\D + 2 \lambda) q(\D) + c\lambda^3 \D^N$. However, a
$2$-cocycle is trivial if and only if it is of the form $(\D + 2
\lambda) q(\D)$, whence the claim.
\end{proof}
\begin{rem}
All of the above non-trivial central extensions of $\Vir$ are
equivalent to one of the form $[L_\lambda L] = (\D + 2 \lambda) L +
\lambda^3 c$, where $\D c = 0$.
\end{rem}

\begin{thm}\lbb{irrcentrext}
A finite non-zero irreducible central extension of $\Vir$ is
isomorphic to that given in \eqref{irreducible}.
\end{thm}
\begin{proof}
We already know that a central extension of $\Vir$ by the
$\CD$-module $C$ is only possible if $C$ is torsion. A torsion
finitely generated $\CD$-module is a finite-dimensional vector
space, on which $\D$ acts as a $\C$-linear endomorphism. Then $C$
decomposes into a direct sum of a submodule on which $\D$ acts
invertibly, and of summands as in Proposition \ref{jordanblock}.

Irreducibility and Proposition \ref{trivial} prove that the summand
on which $\D$ acts invertibly is trivial. On the other hand,
Proposition \ref{jordanblock} shows that we may choose a lifting $L$
of the standard Virasoro generator so that:
$$[L_\lambda L] = (\D + 2 \lambda) L + \lambda^3 c,$$
for some $c \in C$ such that $\D c = 0$. Using again irreducibility
gives $C = \C c$.
\end{proof}

\section{Finite vertex algebras}
\subsection{Commutativity of finite simple vertex algebras}\lbb{wick}

We have seen that coefficients of vertex operators in a vertex
algebra:
\begin{equation}
Y(a,z) = \sum_{j\in \Z} a_{(j)}\,z^{-j-1}
\end{equation}
satisfy the Lie bracket \eqref{liebracket}:
\begin{equation*}
[a_{(m)}, b_{(n)}] = \sum_{j\in \N} {m\choose j}
(a_{(j)}b)_{(m+n-j)}
\end{equation*}
for every $a,b\in V$, $m, n\in \Z$. Multiplying both sides of
\eqref{liebracket} by $\lambda^m z^{-n-1}/m!$ and adding up over all
$m\in \N, n\in \Z$, after applying both sides to $c\in V$, gives
\begin{equation}\lbb{genwickeq}
[a_\lambda Y(b, z) c] = e^{\lambda z} Y([a_\lambda b], z) c + Y(b,
z) [a_\lambda c],
\end{equation}
for all $a,b,c\in V$.

This allows one to explicitly write down the $\lambda$-bracket of a
vertex operator with the normally ordered product of two others --
it suffices to take the constant term in $z$ in both sides -- but
the formula is definitely more useful in the above form. Equation
\eqref{genwickeq} can be used in order to prove the following
statements (see \cite{D}):
\begin{lemma}\lbb{simplelemma}
Let $V$ be a vertex algebra, $U \subset V$ a subspace. Then $[U, V]$
is a (vertex) ideal of $V$.
\end{lemma}

\begin{rem}
It is important to realize that, by Lemma \ref{simplelemma},
elements of the descending sequence \eqref{centralseries} are indeed
ideals of the vertex algebra $V$ and not only of the Lie conformal
algebra $V^\conf$.
\end{rem}

The lemma above has the following immediate and striking
consequence.

\begin{thm}\lbb{simple}
Let $V$ be a non-commutative simple vertex algebra. Then $V^\conf$
is an irreducible central extension of a strongly simple Lie
conformal algebra.
\end{thm}
\begin{rem}
The strong simplicity property is clearly expressed in the proof but
not explicitly stated in \cite{D}.
\end{rem}

We will use Theorem \ref{simple}, and our knowledge of finite simple
Lie conformal algebras, in order to show that all simple vertex
algebra structures over finitely generated $\C[\D]$-modules are
commutative.

\begin{prop}\lbb{novir}
There is no vertex algebra $V$ such that $V = \C[\D]L + \C \1$,
where $L \notin \Tor V$ and $[L_\lambda L] = (\partial + 2 \lambda)
L + c\lambda^3 \1$, for some $c \in \C$.
\end{prop}
\begin{proof}
We proceed by contradiction. We know that
$$[L_\lambda L] = (\partial  + 2 \lambda)L + c\lambda^3 \1,\qquad
[L_\lambda \1] = 0,\qquad [\1_\lambda \1] = 0,$$ and that
$$Y(\1, z)\1 = \1,\qquad Y(\1, z) L = L,\qquad Y(L, z) \1 = e^{z\partial} L.$$
All that we need to determine is $Y(L, z) L$. Let us write
\begin{equation}\lbb{lzl}
Y(L, z) L = a(\partial, z)L + b(z) \1.
\end{equation}
Then \eqref{genwickeq} gives
\begin{align}\lbb{vir}
[L_\lambda Y(L, z) L] & = e^{\lambda z} Y([L_\lambda L], z) L + Y(L,
z) [L_\lambda L],
\end{align}
which, after expanding and comparing coefficients of $L$, yields
\begin{equation}\lbb{virL}
\begin{split}
\left(e^{\lambda z} - 1\right) \frac{d a(\D, z)}{dz} + \left(
2\lambda (e^{\lambda
z} + 1) + \D \right) & \,a(\D, z) =\\
 - c \,\lambda^3 \left(e^{\lambda z} + e^{z\D}\right) & + (\D +
2\lambda) a(\D + \lambda, z).
\end{split}
\end{equation}
Using \eqref{lzl} and substituting $\lambda = -\D/2$, this becomes
\begin{equation}\lbb{diffeq}
\left(e^{-z\D/2} -1\right)\frac{d a(\D,z)}{dz} - \D
e^{-z\D/2}a(\D,z) = \frac{c\D^3}{8}\left(e^{z\D}+e^{-z\D/2}\right),
\end{equation}
i.e., a linear differential equation in $a(\D, z)$ whose solutions
are of the form
$$a(\D,z) = \frac{K(\D)e^{z\D}}{(e^{z\D/2}-1)^2}-\frac{c\D^2}{8}(1 + e^{z\D}).$$
In order for $a(\D,z)$ to be compatible with $[L\,_\lambda L] =
(\partial + 2 \lambda) L$, one needs
$$a(\D,z) = 2/z^2 + \D/z + \mbox{ (regular in $z$)}.$$ This forces
$K(\D) = \D^2/2$, hence the only solution of \eqref{diffeq}
satisfying this additional condition is
\begin{equation}\lbb{adz}
a(\D,z) = \frac{\D^2 e^{z\D}}{2(e^{z\D/2} - 1)^2}-\frac{c\D^2}{8}(1
+ e^{z\D}).
\end{equation}
Checking that this value of $a(\D,z)$ is not a solution of
(\ref{virL}) is a rather lengthy but straightforward computation.
One may also observe that substituting \eqref{adz} into the
right-hand side of \eqref{virL} gives a denominator of the form
$(e^{z(\D + \lambda)/2} -1)^2$ which cannot be obtained from the
left-hand side.
\end{proof}

\begin{thm}\lbb{simpleiscommutative}
Every simple finite vertex algebra is commutative
\end{thm}
\begin{proof}
Let $V$ be a finite simple vertex algebra. By Theorem \ref{simple},
either $V$ is commutative or $V^\conf$ is an irreducible central
extension of a strongly simple Lie conformal algebra. It is then
enough to address the latter case, showing it leads to a
contradiction.

We have seen that every finite strongly simple Lie conformal algebra
is isomorphic to $\Vir$. Moreover, Theorem \ref{irrcentrext} gives a
description of all finite non-zero irreducible central extensions of
$\Vir$. Thus we know that $V = \CD L + \C \1$, with $[L_\lambda L] =
(\partial + 2 \lambda) L + c \lambda^3 \1,$ for some $0 \neq c \in
\C$. Then Proposition \ref{novir} leads to a contradiction.
\end{proof}

The following claim is a technical statement that we will use later
on.

\begin{lemma}\lbb{idealisvertex}
Let $V$ be a vertex algebra, and $M \subset V$ a minimal ideal such
that $M = \C[\D] L + \C k$, where $L$ is a non-torsion element, $\D
k = 0$ and
$$[L_\lambda L] = (\D + 2\lambda)L + \lambda^3 k.$$ Then $M$
can be endowed with a vertex algebra structure by choosing the
vacuum element $\1_M$ to be a suitable scalar multiple of $k$.
\end{lemma}
\begin{proof}
The vacuum element $\1$ of $V$ lies outside of $M$, so the only
thing we need to prove is that we may choose an element inside $M$
whose quantum field acts as the identity on $M$. As $\D k = 0$, then
$Y(k, z)$ does not depend on $z$. Moreover $Y(k, z)k \in M$ is a
torsion element as
$$\D(Y(k, z)k) = Y(\D k, z)k + Y(k, z)(\D k) = 0.$$
Then $Y(k,z)k = \alpha k$ for some $\alpha \in \C$, hence $Y(k -
\alpha\1, z)k = 0$. The element $c = k - \alpha\1$ satisfies
$\partial c = 0$, hence is a torsion element, contained in the
centre of $V$. By Lemma \ref{kercisideal}, $\ker c_{(-1)}$ is an
ideal of $V$ containing $k$, therefore it must contain all of $M$ by
the minimality assumption. Thus, $Y(k - \alpha\1, z)$ has zero
restriction on all of $M$, and $Y(k, z)|_M = \alpha \id_M$. If
$\alpha \neq 0$, then we are done by setting $\1_M = \alpha^{-1} k$.


The case $\alpha = 0$ can be ruled out as follows: we know that
$k_{(-1)}V \subset M$ as $M$ is an ideal containing $k$. Moreover,
we just showed that $k_{(-1)}M = 0$. Now, by \eqref{liebracket}, for
every choice of $a, b \in V$, one has:
\begin{align}
(k_{(-1)}a)_{(m)} (k_{(-1)}b) = k_{(-1)}((k_{(-1)}a)_{(m)} b) +
\sum_{j \geq 0} {m \choose j} ((k_{(-1)}a)_{(j)} k)_{(m+j-1)} b,
\end{align}
where both summands on the right hand side vanish. This shows that
$k_{(-1)} V$ is a subalgebra contained in $M$ in which all products
vanish, therefore $k_{(-1)} V \subset \C k$. In other words, $\C k$
is an ideal of $V$, which contradicts the minimality of $M$.
\end{proof}

\subsection{Solvability of finite vertex algebras}\lbb{solvability}

\begin{lemma}\lbb{nocotorsion}
Let $V$ be a finite vertex algebra, and $S$ be the intersection of
all vertex subalgebras $U \subset V$ such that $\rk U = \rk V$. Then
$V$ is solvable if and only if $S$ is.
\end{lemma}
\begin{proof}
Let $U$ be a vertex subalgebra of $V$ such that $\rk U = \rk V$. $U$
is clearly a subalgebra of $V^\conf$; hence, by Lemma
\ref{cotorsion}, an ideal containing the derived subalgebra of
$V^\conf$.

The intersection $S$ of all such vertex subalgebras is then itself
an ideal of $V^\conf$ containing its derived subalgebra, hence
$V^\conf/S$ is abelian. Therefore, $V$ is solvable if and only if
$S$ is.
\end{proof}

\begin{rem}\lbb{reduction}
Observe that if the vertex subalgebra $S$ in the above lemma is such
that $\rk S = \rk V$, then it is the minimal vertex subalgebra of
$V$ of rank equal to $\rk V$. In particular, $S$ possesses no proper
vertex subalgebras of equal rank.
\end{rem}

\begin{lemma}\lbb{notorsion}
Let $V$ be a finite vertex algebra, and $N$ be the sum of all vertex
ideals of $V$ contained in $\Tor V$. Then $V/N$ has no nonzero
torsion ideal and $V$ is solvable if and only if $V/N$ is. Moreover,
$V$ contains proper vertex subalgebras of rank $\rk V$ if and only
if $V/N$ does.
\end{lemma}
\begin{proof}
The sum of ideals in a vertex algebra is again an ideal. Also, the
sum of torsion elements lies in $\Tor V$. Hence, the sum $N$ of all
vertex ideals of $V$ contained in $\Tor V$ is the maximal such ideal
of $V$. As $\Tor V$ lies in the centre of $V^\conf$, $N^\conf$ is
abelian, so $V$ is solvable if and only if $V/N$ is.

The other claims follow from the correspondence between ideals
(resp. subalgebras) of $V/N$ and ideals (resp. subalgebras) of $V$
containing $N$, and the fact that torsion modules are of zero rank.
\end{proof}

\begin{thm}\lbb{finiteissolvable}
Every finite vertex algebra is solvable.
\end{thm}
\begin{proof}
Assume by way of contradiction that $V$ is a counter-example of
minimal rank. By Lemmas \ref{nocotorsion} and \ref{notorsion} and
Remark \ref{reduction}, we may assume that $V$ has no proper vertex
subalgebra of equal rank, and no non-zero torsion ideal.

Now, observe that if $I \subset V$ is a non-zero ideal with $\rk I <
\rk V$, then $\rk V/I < \rk V$ as $I$ cannot lie in $\Tor V$. Then
the vertex algebras $V/I$ and $I + \C \1$ are both solvable by the
minimality assumption, hence $I^\conf \subset (I + \C \1)^\conf$ is
solvable and $V^\conf$ is an extension of solvable Lie conformal
algebras, a contradiction. Therefore, all non-zero vertex ideals of
$V$ are of the same rank as $V$.

As a consequence, either $V$ is simple, or has a unique non-zero
proper ideal $M$ which is a complement to $\C \1$. Indeed, if $M$ is
a non-zero ideal, then $\rk M = \rk V$, and $M + \C \1$ is a vertex
subalgebra of $V$, hence $M + \C \1 = V$.

We already know that finite simple vertex algebras are commutative,
hence solvable, so it is enough to address the non-simple case. Let
$V$ be a non-solvable finite vertex algebra whose only non-zero
vertex ideal $M \neq V$ is such that $V= M + \C \1$. Then $[V, V] =
[M, M]$ is a vertex ideal of $V$, hence it equals $M$. If $U \subset
M$ is a subspace, so Lemma \ref{simplelemma} shows that $[U, V]$ is
a (vertex) ideal of $V$, hence it equals either $0$ or $M$. As a
consequence, if $u \in M$ then either $u$ is central in $V$ or $[u,
M] = [u, V] = [\C u, V] = M$. This shows that $M^\conf$ is either
strongly simple or a central extension of a strongly simple
conformal algebra. As $[M, M] = M$, the central extension must be
irreducible.

If $M^\conf$ is strongly simple, then we conclude that $V$ is as in
Proposition \ref{novir}, with $c = 0$, hence a contradiction. If, on
the other hand, $M^\conf$ is a central extension of a strongly
simple Lie conformal algebra, then Lemma \ref{idealisvertex} shows
that $M$ can be given a vertex algebra structure which contradicts
Proposition \ref{novir}.
\end{proof}

\section{Conformal adjoint decomposition}\lbb{adjointdecomposition}

\subsection{Finite modules over finite solvable Lie conformal algebras}
In this paper, we will need some basic results from representation
theory of solvable and nilpotent Lie conformal algebras. A
representation of a Lie conformal algebra $R$ is a $\CD$-module $V$
along with a $\lambda$-action $R \otimes V \ni r\otimes v \to
r_\lambda v \in V[\lambda]$ such that

\begin{align}
(\D r)_\lambda v = -\lambda r_\lambda v, & \qquad r_\lambda (\D v) =
(\D + \lambda) r_\lambda v,\\
r_\lambda (s_\mu v) - & s_\mu (r_\lambda v)  = [r_\lambda
s]_{\lambda + \mu} v,
\end{align}
for all $r, s \in R, v \in V$. The action of $r \in R$ on $V$ is
nilpotent if $$r_{\lambda_1} (r_{\lambda_2} ( \dots (r_{\lambda_n}
v) \dots )) = 0$$ for sufficiently large $n$. The following
conformal versions of Engel's and Lie's Theorems were proved in
\cite{DK}.

\begin{thm}\lbb{engel}
Let $R$ be a finite Lie conformal algebra for which every element
$r\in R$ has a nilpotent adjoint action. Then $R$ is a nilpotent Lie
conformal algebra.
\end{thm}

\begin{thm}\lbb{lie}
Let $R$ be a finite Lie solvable conformal algebra, $V$ its finite
module. Then there exists $0 \neq v \in V$ and $\phi: R \ni r \to
\phi_r(\lambda) \in \C[\lambda]$ such that
\begin{equation}
r_\lambda v = \phi_r(\lambda) v,
\end{equation}
for all $r \in R$.
\end{thm}
An element $v$ such as that in Theorem \ref{lie} is a {\em weight
vector}. Then $\phi$ is the {\em weight} of $v$, and it necessarily
satisfies $\phi_{\partial r}(\lambda) = -\lambda \phi_r(\lambda)$.
The set of all weight vectors of a given weight $\phi$, along with
zero, is the {\em weight subspace} $V_\phi$.

\begin{rem}
The statement in \cite{DK} only deals with $R$-modules that are free
as $\CD$-modules, but clearly extends to non-free modules, since
$\Tor V$ is a submodule of $V$ which is killed by $R$.
\end{rem}

\begin{lemma}
Let $V$ be a representation of the Lie conformal algebra $R$. Then
$V_\phi$ is always a vector subspace of $V$. Also, it is a
$\CD$-submodule whenever $\phi \equiv 0$.
\end{lemma}
Now set:
$$V^\phi_{0} = 0, \qquad V^\phi_{i+1} = \left\{ v \in V |\, r_\lambda v - \phi(r) v \in
V^\phi_i \mbox{ for all } r \in R\right\}, i \geq 0.$$ Then
$V^\phi_1 = V_\phi$, and $V^\phi_1 \subset V^\phi_2 \subset ...$ is
an ascending chain of subspaces of $V$. The subspace $\bigcup
V^\phi_i = V^\phi$ is the {\em generalized weight subspace} of
weight $\phi$. Clearly, $r$ acts nilpotently on $V$ exactly when $V$
coincides with the generalized $0$-weight space for the action of
(the Lie conformal algebra generated by) $r$.
\begin{prop}[\cite{DK, BDK}]\lbb{direct}
Let $V$ be a representation of the Lie conformal algebra $R$. Then:
\begin{itemize}
\item $V^\phi$ is a $\CD$-submodule of $V$;
\item $V/V^0$ has no $0$-weight vectors: in particular, it is torsion-free;
\item if $V$ is torsion-free, then $V/V^\phi$ is too;
\item if $\phi \neq \psi$, then $V^\phi \cap V^\psi = 0$;
\item the sum of all generalized weight spaces for the action of $R$
on $V$ is direct.
\end{itemize}
\end{prop}
The sum of all generalized weight spaces may fail to coincide with
the $R$-module $V$. However, the following generalized-weight-space
decomposition, proved in the context of Lie pseudoalgebras
\cite{BDK}, holds for nilpotent Lie conformal algebras.
\begin{thm}\lbb{fitting}
Let $R$ be a finite nilpotent Lie conformal algebra, $V$ its finite
module. Then $V$ decomposes into a direct sum of generalized weight
subspaces for the action of $R$.
\end{thm}
In practice, we will often consider weight spaces and generalized
weight spaces with respect to the action of a single element $s
\in R$.
If $S$ is the subalgebra generated by $s$, we will say a weight for
the action of $S$ on some module $V$ is a {\em weight} of $s$. This
abuse of notation is justified by the fact that in the case $S =
\langle s \rangle = \CD s + S'$, any weight $\phi$ for the action of
$S$ on some module $V$ satisfies $\phi(S') = 0$.

\subsection{Matrix form}\lbb{matrix}

Let $R$ be a Lie conformal algebra, and $V$ be an $R$-module. Then
the map $V \ni v \to r_\lambda v \in V[\lambda]$ is conformal linear
for all $r \in R$. The $\CD$-module structure built on $\Chom(V, V)$
is such that the map $r \mapsto \{ v \mapsto r_\lambda v\}$ is
$\CD$-linear.

One may indeed build up a Lie conformal algebra structure on
$\Chom(V, V)$ in such a way that the above map is always a
homomorphism of Lie conformal algebras. It suffices to define:
\begin{equation}\lbb{commutator}
[f_\lambda g]_\mu v = f_\lambda (g_{\mu-\lambda} v) -
g_{\mu-\lambda} (f_\lambda v),
\end{equation}
whenever $f, g \in \Chom(V, V)$. This Lie conformal algebra
structure is usually denoted by $\gc(V)$, or simply $\gc_n$ when $V$
is a free $\CD$-module of rank $n$. The standard way to represent
elements of $\gc_1$ is by identifying it with $\C[\D, x]$, with the
$\CD$-module structure given via multiplication by $\D$, and the
conformal linear action on $\CD$ given on its free generator $1$ by:
$$x^n_\lambda 1 = (\D + \lambda)^n.$$
Then the $\lambda$-bracket $[p(\D, x)_\lambda q(\D, x)]$ equals
\begin{equation}\lbb{lambdabracket}
p(-\lambda, x + \D + \lambda)q(\D + \lambda, x) - q(\D + \lambda,
x-\lambda)p(-\lambda, x).
\end{equation}

However, in this paper I will employ a different choice, and denote
elements of $\gc_1$ by the effect they have on the free generator.
This identifies $\gc_1$ with $\C[\D, \lambda]$, and has two major
inefficiencies: first of all, the $\CD$-module structure is obtained
via multiplication by $-\lambda$; moreover, as we are already
employing $\lambda$ to denote elements, we will have to compute
$\alpha$- rather than $\lambda$-bracket of elements. However, this
choice is by far more readable than the standard one. The bracket
expressed in \eqref{lambdabracket} then becomes:
\begin{equation}\lbb{alphabracket}
[a(\D, \lambda)_\alpha b(\D, \lambda)] = a(\D, \alpha) b(\D +
\alpha, \lambda-\alpha) - b(\D, \lambda - \alpha) a(\D + \lambda -
\alpha, \alpha).
\end{equation}

Now let $V,W$ be $\CD$-modules. A conformal linear map $f\in
\Chom(V,W)$ is determined by its values on a set of $\CD$-generators
of $V$. If $V$ and $W$ are free, then, for any given choice of
$\CD$-bases $(v^1, ..., v^m)$, $(w^1, ..., w^n)$ of $V$ and $W$
respectively, we can establish a correspondence between $\Chom(V,
W)$ and $n \times m$ matrices with coefficients in $\C[\D,
\lambda]$, similarly to what done above in the case of $\gc_1$.

In general, $\CD$-modules fail to be free. If $M$ is a finitely
generated $\CD$-module, $M$ can be (non-canonically) decomposed as
the direct sum of a free module and of its torsion submodule $\Tor
M$. By Lemma \ref{torsion}, $f \in \Chom(M,N)$ always maps $\Tor M$
to zero.

Since we need to employ a matrix representation for any conformal
linear map $f\in \Chom(M,N)$ between finitely generated modules that
may (and typically will) fail to be free, we can proceed as follows.
Decompose $M$ and $N$ as a direct sum of a free module and their
torsion submodule. If we pick a free $\CD$-basis of the free part,
and a $\C$-basis of the torsion part, we can use this set of
generators to represent conformal linear maps through matrices: we
will call such a set of generators a {\em base}. As a conformal
linear map in $\Chom(M,N)$ always factors via $M/\Tor M$, which is
free, special care is only needed for the treatment of torsion in
the range module.

Note that if we agree that the $\CD$-linear combination expressing
elements of $N$ in terms of a given base is such that coefficients
multiplying torsion elements lie in $\C$ (rather than in $\CD$) then
all coefficients are uniquely determined. This unique expression
enables us to write down a well-behaved matrix representing the
conformal linear map. Matrix coefficients corresponding to torsion
elements then lie in $\C[\lambda]$ rather than in $\C[\D,\lambda]$.

Note that, if $f,g\in \gc(M,M)$ and the matrices representing them
are given by $$F=(f_{ij}(\D,\lambda)), G=(g_{ij}(\D,\lambda)),$$
respectively, then by \eqref{commutator}, the matrix representing
$[f\,_\alpha\,g]$ is given by
\begin{equation}\lbb{matrixcommutator}
F(\D, \alpha)G(\D+\alpha, \lambda - \alpha) - G(\D, \lambda -
\alpha) F(\D + \lambda - \alpha, \alpha),
\end{equation}
where multiplication of matrices is the usual row-by-column product.
Notice that, according to such a matrix representation of conformal
linear maps, Theorem \ref{lie} guarantees the existence of a base in
which matrices representing the action of the solvable Lie conformal
algebra $R$ are simultaneously upper triangular. Similarly, Theorem
\ref{fitting} means that matrices can be put in block diagonal form,
where each block represents the action on a single generalized
weight submodule.

Later, we will call the diagonal entries of a triangular matrix
representing the action of some $s \in S$, $S$ solvable, {\em
eigenvalues} of the element $s$.

\subsection{Adjoint action on a vertex algebra of a Lie conformal
subalgebra}\lbb{decomposition}

In what follows $V$ will be a finite vertex algebra, unless
otherwise stated. If $S$ is a Lie conformal subalgebra of $V$, then
$S$ is solvable by Theorem \ref{finiteissolvable}.
Using formula (\ref{genwickeq}) I want to show that
\begin{thm}\lbb{vertexideal}
If $\psi$ is a non-zero weight for the adjoint action of a Lie
conformal subalgebra $S$ on the finite vertex algebra $V$, then the
generalized weight space $V^\psi$ is a vertex ideal of $V$, and it
satisfies $V^\psi \cdot V^\psi = 0$.
\end{thm}

I will divide the proof of Theorem\,\ref{vertexideal} in a few easy
steps. Let $\alpha$ be a weight for the action of $S$ on $V$,
$\beta$ for its action on $V/V^\alpha$. Denote by $V^{\alpha,
\beta}$ the preimage of $(V/V^\alpha)^\beta$ via the canonical
projection $\pi: V \to V/V^\alpha$. Then we have:
\begin{lemma}\lbb{weightsum}
Let $U\subset V$ be a proper $S$-submodule of $V$ with the property
that $U \cdot V^\psi \subset V^\psi$, and choose an element $w \in
V$ such that $\overline w = [w] \in V/U$ is a weight vector of
weight $\phi$. Then $w \cdot V^\psi \subset V^{\psi, \psi+\phi}$.
\end{lemma}
\begin{proof}
As $[s_\lambda w] = \phi_s(\lambda) w \mod U$, then we have
$s_{(h)}w = \phi_s^h w + u_s^h$, for some $u_s^h \in U$, where the
$\phi_s^h$ are such that $$\phi_s(\lambda) = \sum_h \phi_s^h
\frac{\lambda^h}{h!}.$$ I will prove that
$$w \cdot V^\psi_n \subset V^{\psi,\psi+\phi}$$ by induction on
$n$ -- the basis of induction $n=0$ being trivial, as $V^\psi_0 =
0$.

Let $b\in V^\psi_{n+1}$, and set $s_{(h)}b = \psi_s^h b + v_s^h$
with $v_s^h \in V^\psi_n$. We know that $w_{(N)}b= 0 $ for
sufficiently large $N$. So if $Y(w, z) b\notin V^{\psi,
\psi+\phi}[[z, z^{-1}]]$ we choose $k$ maximal with respect to the
property that $w_{(k)}b \notin V^{\psi,\psi+\phi}$. Let us compute
by means of \eqref{liebracket}:
\begin{equation}
\begin{split}
s_{(m)}(w_{(k)}b) - & w_{(k)}(s_{(m)}b) = \sum_{j=0}^m {m \choose j}
(s_{(j)}w)_{(m+k-j)}b\\ & = (s_{(m)}w)_{(k)}b + \sum_{j=0}^{m-1} {m
\choose j} (s_{(j)}w)_{(m+k-j)}b,
\end{split}
\end{equation}
hence
\begin{equation}
\begin{split}
s_{(m)}(w_{(k)}b) - (\psi_s^m & + \phi_s^m) w_{(k)}b = w_{(k)}v_s^m
+ (u_s^m)_{(k)}b \\ + & \sum_{j=0}^{m-1} {m \choose j}\left(
\phi_s^j (w_{(m+k-j)}b) + (u_s^j)_{(m+k-j)}b\right).
\end{split}
\end{equation}
Now, $v_s^h \in V^\psi_n$, so $Y(w, z)v_s^h \in V^{\psi,
\psi+\phi}[[z, z^{-1}]]$. Also, $u_s^h \in U$, hence $Y(u_s^h, z) b
\in V^\psi[[z, z^{-1}]]$. Moreover, each $w_{(m+k-j)} b$ in the
summation lies in $V^{\psi,\psi+\phi}$ by the maximality of $k$.
Therefore, $(s_{(m)} - (\psi+\phi)_s^m) (w_{(k)}b) \in V^{\psi,
\psi+\phi}$, showing $w_{(k)}b \in V^{\psi, \psi+\phi}$, a
contradiction.
\end{proof}

\begin{lemma}\lbb{weightalone}
Under the same hypotheses as in Lemma \ref{weightsum}, $w \cdot
V^\psi \subset V^\psi$.
\end{lemma}
\begin{proof}
The statement is clear if $\phi = 0$, as $V^{\psi,\psi} = V^\psi$.
Otherwise, choose a base $\{r^i\}$ of $V^{\psi,\psi+\phi}/V^\psi$ on
which the action of $S$ is triangular, and lift it to
$V^{\psi,\psi+\phi}$. Then, if $b\in V^\psi$, we can express $Y(w,
z) b$ as some (depending on $z$) element from $V^\psi$ plus a
$\C[\D]((z))$-linear combination of elements from this base:
\begin{equation}
Y(w, z) b = v(z)+\sum_i A^i(\D,z) r^i.
\end{equation}
My aim is to show that all $A^i$ are zero. I will prove that it is
so for $b \in V^\psi_k$ by induction on $k$. If not all of the $A^i$
are zero, choose $N$ maximal such that $A^N$ is non-zero. Then
(\ref{genwickeq}) gives
$$[s_\lambda Y(w,z)b] = e^{\lambda z} Y([s_\lambda w],z) b + Y(w, z)
[s_\lambda b],$$ and using triangularity of the action of $s$ on the
chosen base, along with the induction assumption, shows that
\begin{equation}\lbb{trans}
(\phi_s(\lambda)+\psi_s(\lambda))A^N(\D+\lambda,z) = (e^{\lambda
z}\phi_s(\lambda)+\psi_s(\lambda))A^N(\D,z),
\end{equation}
as $[s_\lambda b]- \psi_s(\lambda)b$ lies inside $V^\psi_{k-1}$.

Now, since neither $\phi$ nor $\psi$ is identically zero, there must
be some $s$ such that $\phi_s$ and $\psi_s$ are both non-zero. If
for such an $s$ we get $\phi_s+\psi_s = 0$, then $A^N$ must be zero,
giving a contradiction. If instead $\phi_s+\psi_s \neq 0$, then
\begin{equation}
\Gamma(\lambda,z) =
\frac{e^{\lambda z}
\phi_s(\lambda)+\psi_s(\lambda)}{\phi_s(\lambda)+\psi_s(\lambda)}
\end{equation}
is a non-zero element of $\C(\lambda)[[z]]$ satisfying
\begin{equation}\lbb{exp}
\Gamma(\lambda+\mu,z) = \Gamma(\lambda,z)\Gamma(\mu, z).
\end{equation}
It is then easy to show that $\Gamma$'s constant term as a power
series in $z$ must be one. $\Gamma(\lambda,z)$ is indeed of the
form $e^{\lambda \gamma(z)}$ for some power series $\gamma(z) =
\gamma_1 z + \gamma_2 z^2 + ...$

By comparing coefficients of $z$ and $z^2$ in (\ref{exp}) one
concludes that $\phi_s(\lambda)/(\phi_s(\lambda)+\psi_s(\lambda))=0$
or $1$. But this is only possible if either $\phi_s$ or $\psi_s$ is
zero, contrary to the assumption that they are both non-zero. We
obtain a contradiction, which proves that all $A^i$ vanish.
\end{proof}

\begin{lemma}\lbb{ideal}
$V^\psi$ is an ideal of the vertex algebra $V$.
\end{lemma}
\begin{proof}
Let $U$ be maximal among all $S$-submodules of $V$ such that $U
\cdot V^\psi \subset V^\psi$. If $U \neq V$, choose a weight vector
$w$ in $V/U$. Then $(U+\CD w) \cdot V^\psi \subset V^\psi$ by Lemma
\ref{weightalone}, against the maximality of $U$. Hence, $U$ must
equal $V$, and $V^\psi$ is an ideal.
\end{proof}

\begin{lemma}\lbb{sumoweights}
Let $V$ be a (not necessarily finite) vertex algebra, $V^\phi$ and
$V^\psi$ generalized weight subspaces for the adjoint action of the
conformal subalgebra $S$ of $V$. Then $V^\phi \cdot V^\psi \subset
V^{\phi+\psi}.$
\end{lemma}
\begin{proof}
I will show that $V^\phi_i \cdot V^\psi_j \subset V^{\phi+\psi}$ by
induction on $n=i+j$.

Say $v \in V^\phi_i, w\in V^\psi_j$, $i+j=n+1$. Set $s_{(h)}v =
\phi_s^h v + v_s^h, s_{(h)}w = \psi_s^h w + w_s^h$. Then $v_s^h \in
V^\phi_{i-1}, w_s^h \in V^\psi_{j-1}$. If $Y(v, z) w \notin
V^{\phi+\psi}[[z, z^{-1}]]$, then choose a maximal $k$ with the
property that $v_{(k)}w \notin V^{\phi+\psi}$. Then
\begin{equation*}
s_{(m)}(v_{(k)}w) = v_{(k)}(s_{(m)}w) + (s_{(m)}v)_{(k)}w +
\sum_{j=0}^{m-1} {m \choose j} (s_{(j)}v)_{(m+k-j)}w,
\end{equation*}
whence
\begin{equation}\lbb{smvkw}
\begin{split}
s_{(m)}(v_{(k)}w) - (\phi_s^m + & \psi_s^m)(v_{(k)}w) = v_{(k)}w_s^m
+ (v_s^m)_{(k)}w + \\
& \sum_{j=0}^{m-1} {m \choose j}\left(\phi_s^j (v_{(m+n-j)}w) +
(v_s^j)_{(m+n-j)}w\right).
\end{split}
\end{equation}
The right hand side of \eqref{smvkw} lies in $V^{\phi+\psi}$, hence
$v_{(k)}w$ does too, giving a contradiction.
\end{proof}

\begin{proof}[Proof of Theorem \ref{vertexideal}.]
$V^\psi$ is an ideal by Lemma \ref{ideal}. On the other hand, Lemma
\ref{sumoweights} shows $V^\psi \cdot V^\psi \subset V^{2\psi}$. As
$\psi$ is a non-zero weight, Proposition \ref{direct} gives $V^\psi
\cap V^{2\psi} = 0$, hence $V^\psi \cdot V^\psi = 0$.
\end{proof}

Recall that an ideal $I$ of a vertex algebra $V$ is {\em abelian} if
$I \cdot I= 0$, and that $V$ is reduced if it has no abelian ideals
or equivalently if its nilradical is trivial.

\begin{cor}\lbb{nilpotentaction}
Let $V$ be a finite reduced vertex algebra, $N$ be a nilpotent Lie
conformal subalgebra of $V$. Then the adjoint action of $N$ on $V$
is achieved via nilpotent conformal linear maps.
\end{cor}
\begin{proof}
By Theorem~\ref{vertexideal}, any non-zero weight $\phi$ for the
adjoint action of $N$ on $V$ would give an abelian vertex ideal
$V^\phi$. Since this is not possible, the only weight is $0$. But
$N$ is a nilpotent Lie conformal algebra so, by Theorem
\ref{fitting}, every $N$-module decomposes as a direct sum of
generalized weight spaces, showing $V = V^0$. This means that the
action of $N$ on $V$ is nilpotent.
\end{proof}
Corollary \ref{nilpotentaction} has the following immediate
consequence:
\begin{cor}\lbb{nilpocurrent}
If $V$ is a finite reduced vertex algebra, then $\Conf \G$ can arise
as a subalgebra of $V^\conf$ only when $\G$ is a nilpotent Lie
algebra.
\end{cor}
\begin{proof}
Every element $g\in 1 \otimes \G \subset\Conf\G$ spans an abelian
(hence nilpotent) Lie conformal subalgebra of $V$. By
Corollary~\ref{nilpotentaction}, $g$ must act nilpotently on all of
$V$, and in particular on $\Conf \G$ itself. Then $\G$ is a
finite-dimensional Lie algebra on which every element is
$\ad$-nilpotent, and $\G$ is nilpotent by the usual Engel theorem
for Lie algebras.
\end{proof}
\begin{rem}
We observed in Theorem \ref{finiteissolvable} that every finite
vertex algebra is solvable, hence we knew already that $\Conf \G$
arises as a subalgebra of $V^\conf$ only when $\G$ is solvable.
\end{rem}

\section{Nilpotence of finite reduced vertex algebras}\lbb{coreresults}
The main result of this section is the following
\begin{thm}\lbb{main}
Any finite reduced vertex algebra is nilpotent.
\end{thm}
Before we prove this, we show a stronger result characterizing the
conformal adjoint action of elements from $V$.
\begin{lemma}\lbb{peso}
Let $V$ be a finite vertex algebra, $s \in V$. If the (conformal)
adjoint action of $s$ is not nilpotent on $V$, then there exists a
non-zero $\overline s$ whose adjoint action has a weight vector $w$
of non-zero weight.
\end{lemma}
\begin{proof}
In what follows, by \open action of $s$\close, I will always mean
the conformal adjoint action of $s \in V^\conf$ on $V$. Notice that
the finite vertex algebra $V$ is solvable, hence all subalgebras of
$V^\conf$ are solvable Lie conformal algebras, for which the adjoint
action on $V$ satisfies the conditions of Theorem \ref{lie}. In
particular, the subalgebra $\langle s \rangle \subset V^\conf$
generated by $s$ is solvable, and $s$ acts triangularly in a
suitably chosen base of $V$.

If $s$ has a weight vector of non-zero weight then the statement
holds with $\overline s = s$. We can thus assume, without loss of
generality that the only weight of the adjoint action of $s$ on $V$
is zero. Since $s$ does not act nilpotently, $V^0$ cannot equal the
whole $V$ and Proposition \ref{direct} shows that the action of $s$
on $V/V^0$ only has non-zero weights.

Therefore, let us consider a weight vector $\overline w$ in $V/V^0$
of non-zero weight $\phi(\lambda)$. Without loss of generality, we
can assume that the degree of $\phi$ in $\lambda$ be odd. In fact,
we can always replace $s$ by $\partial s$, which must be non-zero,
otherwise $s$ would be a torsion element, with a trivial adjoint
action. Eigenvalues of $\partial s$ are then obtained  by
multiplying those of $s$ by $-\lambda$, and either $\phi(\lambda)$
or $-\lambda \phi(\lambda)$ is of odd degree. Moreover, the adjoint
action of $s$ is nilpotent if and only if that of $\partial s$ is.

So, let $s$ be an element of $V$ for which there exists an element
$\overline w \in V/V^0$ of non-zero weight $\phi(\lambda)$ having
odd degree $n$ in $\lambda$. My plan is to find an element
$\overline s = s+s', s'\in \langle s \rangle'$, and a lifting $w \in
V$ of $\overline w$ in such a way that $w$ will be a (non-zero)
weight vector in $V$ of weight $\phi$ for $\overline s$.

Choose any lifting $w$ of $\overline w$. The $\CD$-submodule $W$
spanned by $V^0$ together with $w$ is preserved by the action of
$S$. Let us fix a base of $W$ consisting of some $S$-triangular base
$(v_1, v_2, ..., v_h)$ for $V^0$ along with the lifting $w$.
According to the matrix representation introduced in
Section~\ref{matrix}, the action of $s$ on $W$ will be represented
by the following matrix:
\begin{equation}
\begin{pmatrix}
0   & * & \hdotsfor[4]{2}& * & X^1(\lambda)\\
0   & 0 & * & \hdotsfor[4]{1}& * & X^2(\lambda) \\
\vdots  & & \ddots & & \vdots \\
0 & \hdotsfor[4]{3}& 0 & X^n(\lambda)\\
0 & \hdotsfor[4]{3} & 0 &  \phi(\lambda)\\
\end{pmatrix},
\end{equation}
where $\phi(\lambda)$ is the weight of $\overline w$. Notice that
coefficients $X^i$ may depend on $\D$ if they refer to an element of
the $\CD$-basis of the free part of $V^0$, but only depend on
$\lambda$ when they refer to basis elements of $\Tor V^0$.

If all of the $X^i$ are zero, then $w$ is a weight vector, and we
are done. If instead some of the $X^i$ are non-zero, let us choose
$i$ to be maximal with the property that $X^i$ is non-zero. I will
show that I can find an element $\overline s = s + s', s' \in S'$
and a lifting $w'$ of $\overline w$ such that, in the matrix
representation of $s+s'$ with respect to the base $(v_1, ..., v_h,
w')$ of $W$, all $X^j, j\geq i$ vanish. An easy induction will then
prove the statement.

If all $X^j, j>i$ are zero and $X^i \neq 0$, then we can compute
the corresponding matrix coefficient in the commutator $[s_\alpha
s]$. By \eqref{matrixcommutator}, it is given by
$$X^i(\D,\alpha) \phi(\lambda-\alpha)
- X^i(\D,\lambda - \alpha) \phi(\alpha).$$
Let us write
$$\phi(\lambda) = \sum_{i=0}^n \phi_i \lambda^i, \qquad
X^i(\D,\lambda) = \sum_{j=0}^m X_j(\D) \lambda^j.$$

Say $m$ is even, and recall that we chose $n$ to be odd. Then the
$i$.th entry in the last column of the matrix representing the
$\alpha^{m+n}$ coefficient (call it $t_1$) in $[s\,_\alpha\, s]$
equals $-2 X_m(\D)\phi_n$. Hence the matrix representing the element
$s_1 = (-\D)^m t_1 / 2\phi_n$ is upper triangular with zero
eigenvalues and the $i$.th entry in the last column is precisely
$-X_m(\D)\lambda^m$, opposite to the highest degree in $\lambda$ of
$X^i(\D,\lambda)$. Thus, the matrix representing $s + s_1$ has the
same eigenvalues as $s$, all $X^j, j>i$ vanish, and the degree in
$\lambda$ of $X^i$ is lower. Notice that $s_1 \in S'$.

On the other hand, if $m$ is odd, assume $m\neq n$. Then the $i$.th
entry in the last column of the matrix representing the coefficient
(call it $t_2$) multiplying $\alpha^{m+n-1}$ in $[s_\alpha s]$ is
given by
$$2 X_m(\D)\phi_{n-1} - 2 X_{m-1}(\D)\phi_n + (n-m)\lambda
X_m(\D)\phi_n.$$ Then the matrix representing the element $s_2 =
(-\D)^{m-1}t_2/(m-n)\phi_n$ is upper triangular with zero
eigenvalues. Moreover, the $i$.th entry in the last column has the
same top degree term in $\lambda$ as $X^i$, with opposite sign. As
before, the element $s+s_2$ has the same eigenvalues as $s$, all
$X^j, j>i$ vanish, and the degree of $X^i$ is lower. Notice that
$s_2 \in S'$ as well. Finally, when $m=n$, it is enough to replace
$w$ with $w + X_m(\D)v_i/\phi_m$ in order to kill the term of top
degree of $X^i$ in the matrix representation of $s$.

By induction, we can then find an $\overline s = s + s', s'\in S'$
and a lifting $w$ of $\overline w$ in such a way that the
corresponding matrix is upper triangular with the same eigenvalues
as $s$, and all $X^j, j\geq i$ vanish.
\end{proof}
\begin{thm}\lbb{key}
Let $V$ be a finite reduced vertex algebra, $s \in V$. Then the
adjoint action of $s$ on $V$ is nilpotent.
\end{thm}
\begin{proof}
If there is some $s\in S$ for which the adjoint conformal action on
$V$ is not nilpotent, then Lemma \ref{peso} finds an element
$\overline s \in \langle s \rangle$ possessing a weight vector of
non-zero weight $\phi$. Then Theorem \ref{vertexideal} shows that
the generalized weight space $V^\phi$ with respect to the conformal
subalgebra generated by $s$ is a non-zero abelian ideal of the
vertex algebra $V$, which is a contradiction, as $V$ is reduced.
\end{proof}
\begin{proof}[Proof of Theorem \ref{main}.]
Every element $v \in V^\conf$ has a nilpotent adjoint conformal
action. By Theorem \ref{engel}, $V^\conf$ is then a nilpotent Lie
conformal algebra.
\end{proof}

\begin{cor}
Let $V$ be a finite vertex algebra, and define $$V^{[0]} = V,\qquad
V^{[n+1]} = [V, V^{[n]}],\,\, n \geq 0.$$ Then all $V^{[n]}$ are
ideals of $V$, and the descending sequence $$V = V^{[0]} \supset
\ldots \supset V^{[n]} \supset \ldots$$ stabilizes on an ideal
contained in the nilradical of $V$.
\end{cor}
\begin{proof}
The quotient of $V$ by its nilradical $N$ is nilpotent, hence the
above sequence for $V/N$ stabilizes to zero. By lifting it back to
$V$, this only happens if the claim holds.
\end{proof}

In other words, every finite vertex algebra can be expressed as an
extension of a nilpotent vertex algebra by a nil-ideal.


\end{document}